\documentclass[11pt, twoside]{amsart}
\usepackage{amsmath, amsthm, amscd, amsfonts, amssymb, graphicx, color}
\usepackage[bookmarksnumbered, plainpages, backref]{hyperref}

\begin{document}

\title{Fibonacci Numbers, Statistical Convergence and Applications}

\author{Murat Kiri\c{s}ci*, Ali Karaisa}

\address{[Murat Kiri\c{s}ci]Department of Mathematical Education, Hasan Ali Y\"{u}cel Education Faculty,
Istanbul University, Vefa, 34470, Fatih, Istanbul, Turkey \vskip 0.1cm }
\email{mkirisci@hotmail.com, murat.kirisci@istanbul.edu.tr}

\address{[Ali Karaisa]Department of Mathematics-Computer Sciences, Necmettin Erbakan University,
Meram Campus, 42090 Meram, Konya, Turkey \vskip 0.1cm }
\email{alikaraisa@hotmail.com; akaraisa@konya.edu.tr}

\thanks{*Corresponding author.}

\begin{abstract}

The purpose of this paper is twofold. First, the definition of new statistical convergence with Fibonacci sequence is given and some fundamental properties
of statistical convergence are examined. Second,  approximation theory worked as a application of the statistical convergence.
\end{abstract}

\maketitle

\pagestyle{plain} \makeatletter
\theoremstyle{plain}
\newtheorem{thm}{Theorem}[section]
\numberwithin{equation}{section}
\numberwithin{figure}{section}  
\theoremstyle{plain}
\newtheorem{pr}[thm]{Proposition}
\theoremstyle{plain}
\newtheorem{exmp}[thm]{Example}
\theoremstyle{plain}
\newtheorem{cor}[thm]{Corollary} 
\theoremstyle{plain}
\newtheorem{defin}[thm]{Definition}
\theoremstyle{plain}
\newtheorem{lem}[thm]{Lemma} 
\theoremstyle{plain}
\newtheorem{rem}[thm]{Remark}
\numberwithin{equation}{section}

\section{Introduction}

\subsection{Densities and Statistical Convergence}
In the theory of numbers, there are many different definitions of density. It is well known that the most popular
of these definitions is asymptotic density. But, asymptotic density does not exist for all sequences. New densities
have been defined to fill those gaps and to serve different purposes.\\

The asymptotic density is one of the possibilities to measure how large a subset of the set of natural number. We
know intuitively that positive integers are much more than perfect squares. Because, every perfect square is
positive and many other positive integers exist besides. However, the set of positive integers is not in fact larger
than the set of perfect squares: both sets are infinite and countable and can therefore be put in one-to-one correspondence.
Nevertheless if one goes through the natural numbers, the squares become increasingly scarce. It is precisely in this case,
natural density help us and makes this intuition precise.\\

Let $A$ is a subset of positive integer. We consider the interval $[1,n]$ and select an integer in this interval, randomly. Then,
the ratio of the number of elements of $A$ in $[1,n]$ to the total number of elements in $[1,n]$ is belong to $A$, probably. For
$n\rightarrow \infty$, if this probability exists, that is this probability tends to some limit, then this limit is used to as
the asymptotic density of the set $A$. This mentions us that the asymptotic density is a kind of probability of choosing a number
from the set $A$.\\

Now, we give some definitions and properties of asymptotic density:\\

The set of positive integers will be denoted by $\mathbb{Z^{+}}$. Let $A$ and $B$ be subsets of $\mathbb{Z}^{+}$. If the symmetric
difference $A\Delta B$ is finite, then, we can say $A$ is asymptotically equal to $B$ and denote $A\sim B$. Freedman and Sember
have introduced the concept of a lower asymptotic density and defined a concept of convergence in density, in \cite{FreeSem}.
\begin{defin}\cite{FreeSem}
Let $f$ be a function which defined for all sets of natural numbers and take values in the interval $[0,1]$.
Then, the function $f$ is said to a lower asymptotic density, if the following conditions hold:
\begin{itemize}
\item[i.] $f(A)=f(B)$, if $A\sim B$,
\item[ii.] $f(A)+f(B)\leq f(A\cup B)$, if $A\cap B=\emptyset$,
\item[iii .] $f(A)+f(B)\leq 1+ f(A\cap B)$, for all $A$,
  \item[iv.]  $f(\mathbb{Z^{+}})=1$.
\end{itemize}
\end{defin}

We can define the upper density based on the definition of lower density as follows:\\

Let $f$ be any density. Then, for any set of natural numbers $A$, the function $\overline{f}$ is said to upper density associated with $f$, if
$\overline{f}(A)=1-f(\mathbb{Z}^{+} \backslash A)$.\\

Consider the set $A\subset \mathbb{Z}^{+}$. If $f(A)=\overline{f}(A)$, then we can say that the set $A$ has natural density
with respect to $f$. The term asymptotic density is often used for the function
\begin{eqnarray*}
d(A)=\liminf_{n\rightarrow\infty}\frac{A(n)}{n},
\end{eqnarray*}
where $A\subset \mathbb{N}$ and $A(n)=\sum_{a\leq n, a\in A}1$. Also the natural density of $A$ is given by $d(A)=\lim_{n}n^{-1}|A(n)|$ where
$|A(n)|$  denotes the number of elements in $A(n)$.\\

The study of statistical convergence was initiated by Fast\cite{Fast}. Schoenberg \cite{Sch} studied statistical convergence as
 a summability method and listed some of the elementary properties of statistical convergence. Both of these mathematicians
  mentioned that if a bounded sequence is statistically convergent to $L$ then it is Ces\`{a}ro summable to $L$. Statistical convergence
  also arises as an example of "convergence in density" as introduced by Buck\cite{Buck}. In \cite{Zygm}, Zygmund called this concept
  "almost convergence" and established relation between statistical convergence and strong summability. The idea of statistical convergence
  have been studied in different branches of mathematics such as number theory\cite{ErTe}, trigonometric series\cite{Zygm}, summability theory\cite{FreeSem},
  measure theory\cite{Mil}, Hausdorff locally convex topological vector spaces\cite{Maddox}. The concept of $\alpha\beta-$statistical convergence was introduced and studied by Aktu\v{g}lu\cite{Aktuglu}. In \cite{VatanAli}, Karakaya and Karaisa have been extended the concept of $\alpha\beta-$statistical convergence. Also,  they have been introduced the concept of weighted $\alpha\beta-$statistical convergence of order $\gamma$, weighted $\alpha\beta-$summmability of order of $\gamma$ and strongly weighted $\alpha\beta-$summable sequences of order $\gamma$, in \cite{VatanAli}. \\

\begin{defin}
A real numbers sequence $x=(x_{k})$ is statistically convergent to $L$ provided that for every $\varepsilon >0$ the set
$\{n\in\mathbb{N}: |x_{n}-L|\geq \varepsilon\}$ has natural density zero. The set of all statistical convergent sequence is denoted by $S$.
In this case, we write $S-\lim x=L$ or $x_{k}\rightarrow L(S)$.\\
\end{defin}

\begin{defin}\cite{Fri}
The sequence $x=(x_{k})$ is statistically Cauchy sequence if for every  $\varepsilon >0$ there is positive integer $N=N(\varepsilon)$ such
that
\begin{eqnarray*}
d\left(\{n\in\mathbb{N}: |x_{n}-x_{N(\varepsilon)}|\geq \varepsilon\}\right)=0.
\end{eqnarray*}

\end{defin}

It can be seen from the definition that statistical convergence is a generalization of the usual of notion of convergence that parallels
the usual theory of convergence.\\

Fridy\cite{Fri} introduce the new notation for facilitate: If $x=(x_{n})$ is a sequence that satisfies some property $P$
for all $n$ except a set of natural density zero, then we say that $x=(x_{n})$ satisfies $P$ for "almost all $n$" and we abbreviate "a.a.n".
In \cite{Fri}, Fridy proved the following theorem:
\begin{thm}
The following statements are equivalent:
\begin{itemize}
\item[i.] $x$ is statistically convergent sequence,
\item[ii.] $x$ is statistically Cauchy sequence,
\item[iii .] $x$ is sequence for which there is a convergent sequence $y$
such that $x_{n}=y_{n}$ for a.a.n.
\end{itemize}
\end{thm}

\subsection{Fibonacci Numbers and Fibonacci Matrix}

The number in the bottom row are called \emph{Fibonacci numbers}, and the number sequence
\begin{eqnarray*}
1,1,2,3,5,8,13,21,34,55,89,144\ldots
\end{eqnarray*}
is the \emph{Fibonacci sequence}\cite{Koshy}. One of the most known and interesting of number sequences is Fibonacci sequence
and it still continues to be of interest to mathematicians. Because, this sequence is an important and useful
tool to expand the mathematical horizon for many mathematician.\\

\begin{defin}
The Fibonacci numbers are the sequence of numbers $(f_{n})$ for $n=1,2,\ldots$ defined by the linear recurrence equation
\begin{eqnarray*}
f_{n}=f_{n-1}+f_{n-2}  \quad n\geq 2,
\end{eqnarray*}
\end{defin}

From this definition, it means that the first two numbers in Fibonacci sequence are either $1$ and $1$ (or $0$ and $1$)
 depending on the chosen starting point of the sequence and all subsequent number is the sum of the previous two. That is,
 we can choose $f_{1}=f_{2}=1$ or $f_{0}=0$, $f_{1}=1$. \\

Fibonacci sequence was initiated in the book \emph{Liber Abaci} of Fibonacci which was written in 1202. However,
the sequence is based on older history. The sequence had been described earlier as Virahanka numbers in Indian mathematics \cite{GooSu}.
In \emph{Liber Abaci}, the sequence starts with $1$, nowadays the sequence begins either with $f_{0}=0$ or with $f_{1}=1$.\\

Some of the fundamental properties of Fibonacci numbers are given as follows:
\begin{eqnarray*}
&&\lim_{n\rightarrow\infty}\frac{f_{n+1}}{f_{n}}=\frac{1+\sqrt{5}}{2}=\alpha, \quad \textrm{(golden ratio)}\\
&&\sum_{k=0}^{n}f_{k}=f_{n+2}-1  \quad  (n\in\mathbb{N}),\\
&&\sum_{k}\frac{1}{f_{k}} \quad \textrm{converges},\\
&&f_{n-1}f_{n+1}-f_{n}^{2}=(-1)^{n+1} \quad (n\geq 1) \quad \textrm{(Cassini formula)}
\end{eqnarray*}
It can be yields $f_{n-1}^{2}+f_{n}f_{n-1}-f_{n}^{2}=(-1)^{n+1}$, if we can substituting for $f_{n+1}$ in Cassini's formula.\\

Let $f_{n}$ be the $n$th Fibonacci number for every $n\in \mathbb{N}$. Then, we define the infinite matrix $\widehat{F}=(\widehat{f}_{nk})$ \cite{Kara1} by

\begin{eqnarray*}
\widehat{f}_{nk}=\left\{\begin{array}{ccl}
-\frac{f_{n+1}}{f_{n}}&, & (k=n-1)\\
\frac{f_{n}}{f_{n+1}}&, & (k=n)\\
0&, & (0\leq k < n-1 \textrm{or} k>n).
\end{array}\right.
\end{eqnarray*}

The Fibonacci Sequence was firstly used in the Theory of Sequence Spaces by Kara and Ba\c{s}ar{\i}r\cite{KaraBas}. Afterward, Kara\cite{Kara1} defined the Fibonacci
difference matrix $\widehat{F}$ by using the Fibonacci sequence $(f_{n})$ for $n\in \{0,1,\ldots\}$ and introduced the new sequence spaces
related to the matrix domain of $\widehat{F}$. \\

Following the \cite{KaraBas} and \cite{Kara1}, high qualified papers are produced with the Fibonacci matrix by many mathematicians.(\cite{AloMur}, \cite{BasBasKara}, \cite{Can1}, \cite{CanKn}, \cite{CanKara}, \cite{CanKay}, \cite{DemKaraBas}, \cite{KaraBasMur}, \cite{KaraSer}, \cite{KaraHan},  \cite{UcBas}).

\subsection{Approximation Theory}

Korovkin type approximation theorems are practical tools to check
whether a given sequence $(A_{n})_{n\geq1}$ of positive linear
operators on $C[a,b]$ of all continuous functions on the real
interval $[a,b]$ is an approximation process. That is, these
theorems present a variety of test functions which provide that the
approximation property holds on the whole space if it holds for
them. Such a property was determined by Korovkin  \cite{Korovkin} in
1953 for the functions $1$, $x$ and $x^{2}$ in the space $C[a,b]$ as
well as for the functions $1$, $\cos$ and $\sin$ in the space of all
continuous $2\pi-$periodic functions on the real line.

Until Gadjiev and Orhan \cite{GadOr} examine, there aren't any study
related to statistical convergence and approximation theory. In
\cite{GadOr}, it is proved a Korovkin type approximation theorems by
using the idea of statistical convergence. Some of examples of
approximation theory and statistical convergence studies can be seen
in \cite{Aktuglu}, \cite{BeMo}, \cite{EdMoNo}, \cite{EdMuKh},
\cite{VatanAli}, \cite{Mohi}, \cite{MuAlo}, \cite{MuAlo1},
\cite{MuVaErG}. Some of examples of approximation theory and
statistical convergence studies can be seen in \cite{Aktuglu},
\cite{BeMo}, \cite{EdMoNo}, \cite{EdMuKh}, \cite{VatanAli},
\cite{Mohi}, \cite{MuAlo}, \cite{MuAlo1}, \cite{MuVaErG}.


\section{Fibonacci Type Statistical Convergence}

Now, we give the general Fibonacci sequence space $X(\widehat{F})$ as follows\cite{Kara1}, \cite{KaraBas}:
Let $X$ be any sequence space and $k\in \mathbb{N}$. Then,
\begin{eqnarray*}
X(\widehat{F})=\left\{x=(x_{k})\in \omega: \left(\widehat{F}x_{k}\right)\in X\right\}
\end{eqnarray*}

It is clear that if $X$ is a linear space, then, $X(\widehat{F})$ is also a linear space.
Kara proved that if $X$ is a Banach space, then, $X(\widehat{F})$ is also a Banach space with the norm
\begin{eqnarray*}
\|x\|_{X(\widehat{F})}=\|\widehat{F}x\|_{X}.
\end{eqnarray*}
Now, we will give lemma which used in proof of Theorem \ref{thm1}. Proof of this lemma is trivial.\\

\begin{lem}\label{lem1}
If $X\subset Y$, then $X(\widehat{F}) \subset Y(\widehat{F})$
\end{lem}

\begin{thm}\label{thm1}
Consider that $X$ is a Banach space and $A$ is a closed subset of $X$. Then,
$A(\widehat{F})$ is also closed in $X(\widehat{F})$.
\end{thm}

\begin{proof}
Since $A$ is a closed subset of $X$, from Lemma \ref{lem1}, then we can write $A(\widehat{F}) \subset X(\widehat{F})$.
$\overline{A(\widehat{F})}$, $\overline{A}$ denote the closure of $A(\widehat{F})$ and  $A$, respectively. To prove the Theorem, we must show that $\overline{A(\widehat{F})}=\overline{A}(\widehat{F})$.\\

Firstly, we take $x\in \overline{A(\widehat{F})}$. Therefore, from 1.4.6 Theorem of \cite{Kreyszig}, there exists a sequence $(x^{n}) \in A(\widehat{F})$
such that $\|x^{n}-x\|_{\widehat{F}}\rightarrow 0$ in $A(\widehat{F})$, for $n\rightarrow\infty$. Thus, $\|(x^{n}_{k})-(x_{k})\|_{\widehat{F}}\rightarrow 0$ as
$n\rightarrow\infty$ in $x\in A(\widehat{F})$ so hat

\begin{eqnarray*}
\sum_{i=1}^{m}\left|x_{i}^{n}-x_{i}\right|+\left\|\widehat{F}(x_{k}^{n})-\widehat{F}(x_{k})\right\|\rightarrow 0
\end{eqnarray*}
for $n\rightarrow\infty$, in $A$. Therefore, $\widehat{F}x \in \overline{A}$ and so $x\in \overline{A}(\widehat{F})$.\\

Conversely, if we take $x\in \overline{A(\widehat{F})}$, then, $x\in A(\widehat{F})$. We know that $A$ is closed. Then $\overline{A(\widehat{F})}=\overline{A}(\widehat{F})$. Hence, $A(\widehat{F})$ is a closed subset of $X(\widehat{F})$.
\end{proof}

From this theorem, we can give the following corollary:

\begin{cor}
If $X$ is a separable space, then, $X(\widehat{F})$ is also a separable space.
\end{cor}

\begin{defin}
 A sequence $x=(x_{k})$ is said to be Fibonacci statistically
convergence(or $\widehat{F}-$statistically convergence) if there is
a  number $L$ such that for every $\epsilon> 0$ the set
$K_{\epsilon}(\widehat{F}):=\{k\leq
n:|\widehat{F}x_{k}-L|\geq\epsilon\}$ has natural density zero,
i.e., $d(K_{\epsilon}(\widehat{F}))=0$. That is

\begin{eqnarray*}
\lim_{n \to \infty}\frac{1}{n}\left|\{k\leq n:
|\widehat{F}x_{k}-L|\geq \epsilon \}\right|=0.
\end{eqnarray*}

In this case we write $d(\widehat{F})-\lim x_{k}=L$ or $x_{k}\rightarrow L(S(\widehat{F}))$. The set of
$\widehat{F}-$statistically convergent sequences will be denoted by
$S(\widehat{F})$. In the case $L=0$, we will write
$S_{0}(\widehat{F})$.

\end{defin}

\begin{defin}\label{defin2}
Let $x=(x_{k})\in \omega$. The sequence $x$ is said to be $\widehat{F}-$statistically Cauchy if
there exists a number $N=N(\varepsilon)$ such that
\begin{eqnarray*}
\lim_{n \to \infty}\frac{1}{n}\left|\{k\leq n:
|\widehat{F}x_{k}-\widehat{F}x_{N}|\geq \epsilon \}\right|=0
\end{eqnarray*}
for every $\varepsilon >0$.
\end{defin}

\begin{thm}
If $x$ is a $\widehat{F}-$statistically convergent sequence, then $x$ is a $\widehat{F}-$statistically Cauchy sequence.
\end{thm}

\begin{proof}
Let $\varepsilon >0$. Assume that $x_{k}\rightarrow L(S(\widehat{F}))$. Then, $|\widehat{F}x_{k}-L|<\varepsilon / 2$
for almost all $k$. If $N$ is chosen so that $|\widehat{F}x_{N}-L|<\varepsilon / 2$, then we have
$|\widehat{F}x_{k}-\widehat{F}x_{N}|< |\widehat{F}x_{k}-L|+|\widehat{F}x_{N}-L|<\varepsilon / 2 + \varepsilon / 2 =\varepsilon$
for almost all $k$. It means that $x$ is $\widehat{F}-$statistically Cauchy sequence.
\end{proof}

\begin{thm}
If $x$ is a sequence for which there is a $\widehat{F}-$statistically convergent sequence $y$ such that
$\widehat{F}x_{k}=\widehat{F}y_{k}$ for almost all $k$, then $x$ is $\widehat{F}-$ statistically convergent sequence.
\end{thm}

\begin{proof}
Suppose that $\widehat{F}x_{k}=\widehat{F}y_{k}$ for almost all $k$ and $y_{k}\rightarrow L(S(\widehat{F}))$. Then, $\varepsilon >0$
and for each $n$, $\{k \leq n: |\widehat{F}x_{k}-L|\geq \varepsilon\}\subseteq \{k \leq n: \widehat{F}x_{k}\neq \widehat{F}y_{k}\} \cup
\{k \leq n: |\widehat{F}x_{k}-L|\leq \varepsilon\}$. Since $y_{k}\rightarrow L(S(\widehat{F}))$, the latter set contains a fixed number of
integers, say $g=g(\varepsilon)$. Therefore, for $\widehat{F}x_{k}=\widehat{F}y_{k}$ for almost all $k$,
\begin{eqnarray*}
\lim_{n}\frac{1}{n}\left|\{k\leq n: |\widehat{F}x_{k}-L|\geq \varepsilon\}\right|\leq
\lim_{n}\frac{1}{n}\left|\{k\leq n: \widehat{F}x_{k}\neq \widehat{F}y_{k} \}\right|+\lim_{n}\frac{g}{n}=0.
\end{eqnarray*}
Hence $x_{k}\rightarrow L(S(\widehat{F}))$.
\end{proof}

\begin{defin}\label{defin1}\cite{FriOrh}
A sequence $x=(x_{k})$ is said to be statistically bounded if there exists some $L\geq 0$ such that
\begin{eqnarray*}
d\left(\{k: |x_{k}|>L\}\right)=0, \quad \textrm{i.e.,} \quad |x_{k}|\leq L \quad \textrm{a.a.k}
\end{eqnarray*}
\end{defin}

By $m_{0}$, we will denote the linear space of all statistically bounded sequences. Bounded sequences
are obviously statistically bounded as the empty set has zero natural density. However, the converse
is not true. For example, we consider the sequence
\begin{eqnarray*}
x_{n} = \left\{ \begin{array}{ccl}
n&, & (\textrm{k is a square}),\\
0&, & (\textrm{k is not a square}).
\end{array} \right.
\end{eqnarray*}
Clearly $(x_{k})$ is not a bounded sequence. However, $d(\{k: |x_{k}|>1/2\})=0$, as the of squares has zero natural density and hence
$(x_{k})$ is statistically bounded\cite{BharGup}.\\

\begin{pr}\label{prop1}\cite{BharGup}
Every convergent sequence is statistically bounded.
\end{pr}

Although a statistically convergent sequence does not need to be bounded(cf. \cite{BharGup},\cite{FriDemOrh}), the following proposition
shows that every statistical convergent sequence is statistically bounded.

\begin{pr}\label{prop2}\cite{BharGup}
Every statistical convergent sequence is statistically bounded.
\end{pr}

Now, using the Propositions \ref{prop1}, \ref{prop2}, we can give the following corollary:

\begin{cor}
Every $\widehat{F}-$statistical convergent sequence is $\widehat{F}-$statistically bounded.
\end{cor}

Denote the set of all $\widehat{F}-$bounded sequences of real numbers by $m(\widehat{F})$ \cite{Kara1}.
Based on the Definition \ref{defin1} and the descriptions of $m_{0}$ and $m(\widehat{F})$, we can denote the set of all $\widehat{F}-$bounded statistically
convergent sequences of real numbers by $m_{0}(\widehat{F})$.\\

Following theorem can be proved from Theorem 2.1 of \cite{salat} and Theorem \ref{thm1}.

\begin{thm}\label{boun1}
The set of $m_{0}(\widehat{F})$ is a closed linear space of the linear normed space $m(\widehat{F})$.
\end{thm}

\begin{thm}
The set of $m_{0}(\widehat{F})$ is a nowhere dense set in $m(\widehat{F})$.
\end{thm}

\begin{proof}
According to \cite{salat} that every closed linear subspace of an arbitrary linear normed space $E$,
different from $E$ is a nowhere dense set in $E$. Hence, on account of Theorem \ref{boun1}, it suffices
to prove that $m_{0}(\widehat{F})\neq m(\widehat{F})$. But this is evident, consider the sequence
\begin{eqnarray*}
x_{n} = \left\{ \begin{array}{ccl}
1&, & (\textrm{n is odd}),\\
0&, & (\textrm{n is even}).
\end{array} \right.
\end{eqnarray*}
Then, $x\in m(\widehat{F})$, but does not belong to $m_{0}(\widehat{F})$.
\end{proof}

$\omega$ denotes the Fr\'{e}chet metric space of all real sequences with the metric $d_{\omega}$,
\begin{eqnarray*}
d_{\omega}=\sum_{k=1}^{\infty}\frac{1}{2^{k}}\frac{|x_{k}-y_{k}|}{1+|x_{k}-y_{k}|}
\end{eqnarray*}
where $x=(x_{k}), y=(y_{k})\in \omega$ for all $k=1,2,\cdots$.

\begin{thm}
The set of $\widehat{F}-$statistically convergent sequences is dense in the space $\omega$.
\end{thm}

\begin{proof}
If $x=(x_{k})\in S(\widehat{F})$ (for all $k$) and the sequence $y=(y_{k})$ (for all $k$) of real numbers
differs from $x$ only in a finite number of terms, then evidently $y\in S(\widehat{F})$, too. From this statement follows at once on
the basis of the definition of the metric in $\omega$.
\end{proof}

\begin{thm}The following statements are hold.
\begin{itemize}
\item[i.] The inclusion $c(\widehat{F})\subset S(\widehat{F})$ is strict.
\item[ii.] $S(\widehat{F})$ and $\ell_{\infty}(\widehat{F})$, overlap but neither one contains the other.
\item[iii .] $S(\widehat{F})$ and $\ell_{\infty}$, overlap but neither one contains the other.
  \item[iv.]  $S$ and $S(\widehat{F})$, overlap but neither one contains the other.
      \item[v.] $S$ and $c(\widehat{F})$, overlap but neither one contains the other.
     \item[vi.]  $S$ and $c_{0}(\widehat{F})$, overlap but neither one contains the other.
     \item[vii.] $S$ and $\ell_{\infty}(\widehat{F})$, overlap but neither one contains the other.
\end{itemize}
\end{thm}

\begin{proof}
i) Since $c\subset S$, then $c(\widehat{F})\subset S(\widehat{F})$.We choose
\begin{eqnarray}\label{sampleseq}
\widehat{F}x_{n}=(f_{n+1}^{2})=(1,2^{2},3^{2},5^{2},\ldots),
\end{eqnarray}
Since $f_{n+1}^{2}\rightarrow \infty$ as $k\rightarrow\infty$ and $\widehat{F}x=(1,0,0,\ldots)$,
 then, $\widehat{F}x\in S$, but is not in the space $c$, that is, $\widehat{F} \not \in c$.\\

For the other items, firstly, using the inclusion relations in \cite{Kara1}. It is obtained that the inclusions $c\subset S(\widehat{F})$, $c\subset c(\widehat{F})$, $c\subset m(\widehat{F})$, $c\subset S$, $c\subset \ell_{\infty}$ and $c\cap c_{0}(\widehat{F})\neq \phi$ are hold. Then, we  seen that $S(\widehat{F})$ and $\ell_{\infty}(\widehat{F})$, $S(\widehat{F})$ and $\ell_{\infty}$, $S$ and $S(\widehat{F})$, $S$ and $c(\widehat{F})$, $S$ and $c_{0}(\widehat{F})$, $S$ and $\ell_{\infty}(\widehat{F})$ is overlap.\\

ii) We define $\widehat{F}x=\widehat{F}x_{n}$ from (\ref{sampleseq}). Then, $\widehat{F}x\in S$, but $\widehat{F}x$ is not in $\ell_{\infty}$. Now we choose $u=(1,0,1,0,\ldots)$. Then,  $u\in \ell_{\infty}(\widehat{F})$ but $u\not \in S(\widehat{F})$.\\

iii) The proof is the same as (ii).\\

iv) Define
\begin{eqnarray*}
x_{n} = \left\{ \begin{array}{ccl}
1&, & (\textrm{n is a square}),\\
0&, & (\textrm{otherwise}).
\end{array} \right.
\end{eqnarray*}
Then $x\in S$ but $x \in S(\widehat{F})$.
Conversely, if we take $u=(n)$, then $u\not \in S$ but $x \in S(\widehat{F})$.\\

(v), (vi) and (vii) are proven similar to (iv).

\end{proof}

\section{Applications}

\subsection{Approximation by $\widehat{F}-$statistically convergence}

In this section, we get an analogue of classical Korovkin Theorem by
using the concept $\widehat{F}-$statistically convergence.

Let $F(\mathbb{R})$ denote the linear space of real value function
on $\mathbb{R}$. Let $C(\mathbb{R})$ be  space of all real-valued
continuous functions $f$ on $\mathbb{R}$. It is well known that
$C(\mathbb{R})$ is Banach space with the norm given as follows:
\begin{eqnarray*}
\parallel f\parallel_{\infty}=\sup_{x\in \mathbb{R}}|f(x)|,\,\,f\in C(\mathbb{R})
\end{eqnarray*}
and we denote $C_{2\pi}(\mathbb{R})$ the space of all $2\pi-$
periodic functions $f\in C(\mathbb{R})$ which is Banach space with
the norm given by
\begin{eqnarray*}
\parallel f\parallel_{2\pi}=\sup_{t\in \mathbb{R}}|f(t)|,\,\,f\in
C(\mathbb{R}).
\end{eqnarray*}
We say $A$ is positive operator, if for every non-negative $f$ and
$x\in I$, we have $A(f,x)\geq 0$, where $I$ is any given interval on
the real semi-axis. The first and second classical Korovkin
approximation theorems states as follows (see \cite{Gadziev}, \cite{Korovkin})
\begin{thm}
Let $(A_{n})$ be a sequence of positive linear operators from
$C[0,1]$ in to $F[0,1]$. Then
\begin{eqnarray*}
\lim_{n \to \infty}\parallel A_{n}(f,x)-f(x)\parallel_{C[a,b]}=0
\Leftrightarrow \lim_{n \to \infty}\parallel A_{n}(e_{i},
x)-e_{i}\parallel_{C[a,b]}=0,
\end{eqnarray*}
where $e_{i}=x^{i}$, $i=0,1,2$.
\end{thm}
\begin{thm}\label{kor2}
Let $(T_{n})$ be sequence of positive linear operators from
$C_{2\pi}(\mathbb{R})$ into $F(\mathbb{R})$. Then

\begin{eqnarray*}
\lim_{n \to \infty}\parallel T_{n}(f,x)-f(x)\parallel_{2\pi}=0
\Leftrightarrow \lim_{n \to \infty}\parallel T_{n}(f_{i},
x)-f_{i}\parallel_{2\pi}=0,\,\, i=0,1,2
\end{eqnarray*}
where $f_{0}=1$, $f_{1}=sinx$ and $f_{2}=cosx$.

\end{thm}
Our main Korovkin type theorem is given as follows:
\begin{thm}\label{thm}
Let $(L_{k})$ be a sequence of positive linear operator from
$C_{2\pi}(\mathbb{R})$ into $C_{2\pi}(\mathbb{R})$. Then for all
$f\in C_{2\pi}(\mathbb{R})$

\begin{eqnarray}\label{k0}
d(\widehat{F})-\lim_{k \to \infty}\parallel
L_{k}(f,x)-f(x)\parallel_{2\pi}=0
\end{eqnarray}
if and only if
\begin{eqnarray}
&&d(\widehat{F})-\lim_{k \to \infty}\parallel L_{k}(1,x)-1\parallel_{2\pi}=0\label{k1},\\
&&d(\widehat{F})-\lim_{k \to \infty}\parallel
L_{k}(\sin t,x)-\sin x\parallel_{2\pi}=0,\label{k2}\\
&&d(\widehat{F})-\lim_{k \to \infty}\parallel L_{k}(\cos t,x)-\cos
x\parallel_{2\pi}=0.\label{k3}
\end{eqnarray}
\end{thm}
\begin{proof}
As $1,\sin x,\cos x\in C_{2\pi}(\mathbb{R})$, conditions
(\ref{k1})-(\ref{k3}) follow immediately from (\ref{k0}). Let the
conditions (\ref{k1})-(\ref{k3}) hold and $I_{1}=(a,a+2\pi)$ be any
subinterval of length $2\pi$ in $\mathbb{R}$. Let fixed  $x\in
I_{1}$. By the properties of function $f$ , it follows that for
given $\varepsilon> 0$ there exists $\delta=\delta(\epsilon)>0$ such
that
\begin{eqnarray}\label{k4}
|f(x)-f(t)|<\varepsilon,\,\,\textrm{whenever
}\,\,\forall|t-x|<\delta.
\end{eqnarray}
If $|t-x|\geq\delta$, let us assume that $t\in
(x+\delta,2\pi+x+\delta)$. Then we obtain that

\begin{eqnarray}\label{k5}
|f(x)-f(t)|\leq2\parallel f\parallel_{2\pi}\leq\frac{2\parallel
f\parallel_{2\pi}}{\sin^{2}\left(\frac{\delta}{2}\right)}\psi (t),
\end{eqnarray}
where $\psi (t)=\sin^{2}\left(\frac{t-x}{2}\right)$.

By using (\ref{k4}) and (\ref{k5}), we  have
\begin{eqnarray*}
|f(x)-f(t)|<\varepsilon+ \frac{2\parallel
f\parallel_{2\pi}}{\sin^{2}\left(\frac{\delta}{2}\right)}\psi (t).
\end{eqnarray*}
This implies that
\begin{eqnarray*}
-\varepsilon-\frac{2\parallel
f\parallel_{2\pi}}{\sin^{2}\left(\frac{\delta}{2}\right)}\psi
(t)<f(x)-f(t)<\varepsilon+ \frac{2\parallel
f\parallel_{2\pi}}{\sin^{2}\left(\frac{\delta}{2}\right)}\psi (t).
\end{eqnarray*}
By using the positivity and linearity of $\{L_{k}\}$, we get

\begin{eqnarray*}
L_{k}(1,x)\left(-\varepsilon \frac{2\parallel
f\parallel_{2\pi}}{\sin^{2}\left(\frac{\delta}{2}\right)}\psi
(t)\right)<L_{k}(1,x)\left(f(x)-f(t)\right)<
L_{k}(1,x)\left(\varepsilon+ \frac{2\parallel
f\parallel_{2\pi}}{\sin^{2}\left(\frac{\delta}{2}\right)}\psi
(t)\right)
\end{eqnarray*}
where $x$ is fixed  and so $ f(x)$ is constant number. Therefore,
\begin{eqnarray}\label{k6}
-\varepsilon L_{k}(1,x)-\frac{2\parallel
f\parallel_{2\pi}}{\sin^{2}\left(\frac{\delta}{2}\right)}L_{k}(\psi
(t),x)< L_{k}(f,x)-f(x)L_{k}(1,x)
\end{eqnarray}
\begin{eqnarray*}
<\varepsilon
L_{k}(1,x)+\frac{2\parallel
f\parallel_{2\pi}}{\sin^{2}\left(\frac{\delta}{2}\right)}L_{k}(\psi
(t),x).
\end{eqnarray*}
On the other hand, we get
\begin{eqnarray}\label{k7}
L_{k}(f,x)-f(x)&=&L_{k}(f,x)-f(x)L_{k}(1,x)+f(x)L_{k}(1,x)-f(x)\nonumber\\
&=&L_{k}(f,x)-f(x)L_{k}(1,x)-f(x)L_{k}+f(x)[L_{k}(1,x)-1]\label{k7}.
\end{eqnarray}
By inequality (\ref{k6}) and (\ref{k7}), we obtain
\begin{eqnarray}\label{k8}
L_{k}(f,x)-f(x)&<&\varepsilon L_{k}(1,x)+\frac{2\parallel
f\parallel_{2\pi}}{\sin^{2}\left(\frac{\delta}{2}\right)}L_{k}(\psi
(t),x)\\ \nonumber &&+f(x)+f(x)[L_{k}(1,x)-1].
\end{eqnarray}
Now, we compute second moment
\begin{eqnarray*}
L_{k}(\psi(t),x)=L_{k}\left(\sin^{2}\left(\frac{x-t}{2}\right),x\right)&=&L_{k}\left(\frac{1}{2}(1-\cos t \cos x-\sin x \sin t),x\right)\\
&=&\frac{1}{2}[L_{k}(1,x)-\cos xL_{k}(\cos t,x)-\sin xL_{k}(\sin t ,x)]\\
&=&\frac{1}{2}\{L_{k}(1,x)-\cos x[L_{k}(\cos t,x)-\cos x]\\
&-&\sin x[L_{k}(\sin t ,x)-\sin x]\}.
\end{eqnarray*}
By (\ref{k8}), we have
\begin{eqnarray*}\label{k9}
L_{k}(f,x)-f(x)&<&\varepsilon L_{k}(1,x)+\frac{2\parallel
f\parallel_{2\pi}}{\sin^{2}\left(\frac{\delta}{2}\right)}\frac{1}{2}\{L_{k}(1,x)\nonumber\\
&&-\cos x[L_{k}(\cos t,x)-\cos x]-\sin x[L_{k}(\sin t ,x)-\sin x]\}+f(x)(L_{k}(1,x)-1)\nonumber\\
&=&\varepsilon[L_{k}(1,x)-1]+\varepsilon+f(x)(L_{k}(1,x)-1)\nonumber\\
&&+\frac{\parallel
f\parallel_{2\pi}}{\sin^{2}\left(\frac{\delta}{2}\right)}\{L_{k}(1,x)-\cos
x[L_{k}(\cos t,x)-\cos x]\nonumber\\
&-&\sin x[L_{k}(\sin t ,x)-\sin x]\}.
\end{eqnarray*}
So, from above inequality, one can see that
\begin{eqnarray*}
|L_{k}(f,x)-f(x)|&\leq&\varepsilon+\left(\varepsilon+|f(x)|+\frac{\parallel
f\parallel_{2\pi}}{\sin^{2}\left(\frac{\delta}{2}\right)}\right)|L_{k}(1,x)-1|
\\ \nonumber &&+\frac{\parallel
f\parallel_{2\pi}}{\sin^{2}\left(\frac{\delta}{2}\right)}[|\cos
x||L_{k}(\cos t,x)-\cos x|\nonumber\\&&+|\sin x|L_{k}(\sin t,x)-\sin
x|] \leq\varepsilon+\left(\varepsilon+|f(x)|+\frac{\parallel
f\parallel_{2\pi}}{\sin^{2}\left(\frac{\delta}{2}\right)}\right)|L_{k}(1,x)-1|\\
\nonumber&&+\frac{\parallel
f\parallel_{2\pi}}{\sin^{2}\left(\frac{\delta}{2}\right)}[|L_{k}(\cos
t,x)-\cos x|+|\sin x|L_{k}(\sin t ,x)-\sin x|].
\end{eqnarray*}
Because of $\varepsilon$ is arbitrary, we obtain
\begin{eqnarray*}
\parallel L_{k}(f,x)-f(x)\parallel _{2\pi}&\leq&
\varepsilon+R\bigg(\parallel L_{k}(1,x)-1\parallel _{2\pi}+\parallel
L_{k}(\cos t,x)-\cos t\parallel _{2\pi}\\&&+\parallel L_{k}(\sin
t,x)-\sin x\parallel _{2\pi}\bigg)
\end{eqnarray*}
where $R=\max\left(\varepsilon+\parallel
f\parallel_{2\pi}+\frac{\parallel
f\parallel_{2\pi}}{\sin^{2}\left(\frac{\delta}{2}\right)},\frac{\parallel
f\parallel_{2\pi}}{\sin^{2}\left(\frac{\delta}{2}\right)}\right)$.\\
Finally, replacing $L_{k}(.,x)$ by
$T_{k}(.,x)=\widehat{F}L_{k}(.,x)$ and for $\varepsilon^{'}>0$, we
can write
\begin{eqnarray*}
\mathcal{A}:&=&\left\{k \in \mathbb{N}:\parallel
T_{k}(1,x)-1\parallel _{2\pi}+\parallel T_{k}(\sin t,x)-\sin
x\parallel _{2\pi}+\parallel T_{k}(\cos t,x)-\cos x\parallel
_{2\pi}\geq\frac{\varepsilon^{'}}{R}\right\},\\
\mathcal{A}_{1}:&=&\left\{k \in \mathbb{N}:\parallel
T_{k}(1,x)-1\parallel _{2\pi}\geq\frac{\varepsilon^{'}}{3R}\right\},\\
\mathcal{A}_{2}:&=&\left\{k \in \mathbb{N}:\parallel T_{k}(\sin
t,x)-\sin x\parallel _{2\pi}\geq\frac{\varepsilon^{'}}{3R}\right\},\\
\mathcal{A}_{3}:&=&\left\{k \in \mathbb{N}:\parallel T_{k}(\cos
t,x)-\cos x\parallel _{2\pi}\geq\frac{\varepsilon^{'}}{3R}\right\}.
\end{eqnarray*}
Then, $\mathcal{A}\subset \mathcal{A}_{1}\cup \mathcal{A}_{2}\cup
\mathcal{A}_{3}$, so we have $d(\mathcal{A})\leq
d(\mathcal{A}_{1})+d(\mathcal{A}_{2})+d(\mathcal{A}_{3})$. Thus, by
conditions (\ref{k1})-(\ref{k3}), we obtain
\begin{eqnarray*}
d(\widehat{F})-\lim_{k \to \infty}\parallel
L_{k}(f,x)-f(x)\parallel_{2\pi}=0.
\end{eqnarray*}
which completes the proof.
\end{proof}

We remark that our Theorem \ref{thm} is stronger than Theorem
\ref{kor2} as well as Theorem of  Gadjiev and Orhan \cite{GadOr}.
For this purpose, we get the following example:
\begin{exmp}
For $n\in\mathbb{N}$, denote by $S_{n}(f)$ the $n-$ partial sum of
the Fourier series of $f$, that is
\begin{eqnarray*}
S_{n}(f,x)=\frac{1}{2}a_{0}(f)+\sum_{k=0}^{n}a_{k}(f)\cos
kx+b_{k}(f)\sin k x.
\end{eqnarray*}
For $n\in\mathbb{N}$, we get
\begin{eqnarray*}
F_{n}(f,x)=\frac{1}{n+1}\sum_{k=0}^{n}S_{k}(f).
\end{eqnarray*}
A standard calculation  gives that for every $t\in\mathbb{R}$
\begin{eqnarray*}
F_{n}(f,x)=\frac{1}{2\pi}\int_{-\pi}^{\pi}f(t)\varphi_{n}(x-t)dt,
\end{eqnarray*}
where
\begin{eqnarray*}
\varphi_{n}(x) = \left\{ \begin{array}{ll}
\frac{\sin^{2}((n+1)(x-t)/2)}{(n+1)\sin^{2}((x-t)/2)}& \textrm{ if $x$ is not a multiple of $2\pi$ },\\
n+1 & \textrm{ if $x$ is a multiple of $2\pi$ }.
\end{array} \right.
\end{eqnarray*}
The sequence $(\varphi_{n})_{n\in\mathbb{N}}$ is a positive kernel
which is called the Fej$\acute{e}$r kernel, and  corresponding
$F_{n}$ for $n\geq 1$ are called Fej$\acute{e}$r convolution
operators.

We define the sequence of linear operators as
$K_{n}:C_{2\pi}(\mathbb{R})\longrightarrow C_{2\pi}(\mathbb{R})$
with $K_{n}(f,x)=(1+y_{n})F_{n}(f,x)$, where
$y=(y_{n})=(f^{2}_{n+1})$ . Then, $K_{n}(1,x)=1$, $K_{n}(\sin
t,x)=\frac{n}{n+1}\sin x$ and $K_{n}(\cos t,x)=\frac{n}{n+1}\cos x$
and sequence $(K_{n})$ satisfies the conditions
(\ref{k1})-(\ref{k3}). Therefore, we get
\begin{eqnarray*}
 d(\widehat{F})-\lim_{k \to \infty}\parallel
K_{n}(f,x)-f(x)\parallel_{2\pi}=0.
\end{eqnarray*}
On the other hand, one can see that $(K_{n})$ does not satisfy
Theorem \ref{kor2} as well as Theorem of  Gadjiev and Orhan
\cite{GadOr}, since $\widehat{F}y=(1,0,0,\ldots)$, the sequence $y$
is $\widehat{F}-$statistical convergence to $0$. But the sequence
$y$ neither convergent nor statistical convergent.
\end{exmp}

\subsection{Rate of $\widehat{F}-$statistical convergence} In this
section, we estimate rate of $\widehat{F}-$statistical convergence
of a sequence of positive linear operators defined
$C_{2\pi}(\mathbb{R})$ into $C_{2\pi}(\mathbb{R})$. Now, we give
following definition
\begin{defin}\label{defn}
Let $(a_{n})$ be a positive non-increasing sequence. We say that the
sequence $x=(x_{k})$ is $\widehat{F}-$statistical convergence to $L$
with the rate $o(a_{n})$ if for every, $\varepsilon>0$
\begin{eqnarray*}
\lim_{n \to \infty}\frac{1}{u_{n}}\left|\left\{k\leq
n:|\widehat{F}x-\ell|\geq\epsilon\right\}\right|=0.
\end{eqnarray*}
At this stage, we can write $x_{k}-L=d(\widehat{F})-o(u_{n})$.
\end{defin}
As usual we have the following auxiliary result.
\begin{lem}
Let $(a_{n})$ and $(b_{n})$  be two positive non-increasing
sequences. Let $x=(x_{k})$ and $y=(y_{k})$ be two sequences such
that $x_{k}-L_{1}=d(\widehat{F})-o(a_{n})$ and
$y_{k}-L_{1}=d(\widehat{F})-o(b_{n})$. Then we have
\begin{enumerate}
\item[(i)] $\alpha(x_{k}-L_{1})=d(\widehat{F})-o(a_{n})$ for any
scalar $\alpha$,\\
\item[(ii)] $(x_{k}-L_{1})\pm (y_{k}-L_{2})=d(\widehat{F})-o(c_{n})$,\\
\item[(iii)] $(x_{k}-L_{1})(y_{k}-L_{2})=d(\widehat{F})-o(a_{n}b_{n}),$
\end{enumerate}
where $c_{n}= \max\{a_{n},b_{n}\}$.
\end{lem}

For $\delta> 0$, the modulus of continuity of $f$,
$\omega(f,\delta)$ is defined by
\begin{equation*}
\omega(f,\delta)=\sup_{|x-y|<\delta }|f(x)-f(y)|.
\end{equation*}
It is well-known that for a function $f \in C[a,b]$,
\begin{equation*}
\lim_{n\to 0^{+}}\omega(f,\delta)=0
\end{equation*}
for any $\delta>0$
\begin{equation}\label{r1}
|f(x)-f(y)|\leq \omega(f,\delta)\left(\frac{|x-y|}{\delta}+1\right).
\end{equation}

\begin{thm}
Let $(L_{k})$ be sequence of positive linear operator from
$C_{2\pi}(\mathbb{R})$ into $C_{2\pi}(\mathbb{R})$. Assume that
\begin{eqnarray*}
(i)&&\ \parallel L_{k}(1,x)-x\parallel
_{2\pi}=d(\widehat{F})-o(u_{n}),\\
(ii)&&\omega(f,\theta_{k})=d(\widehat{F})-o(v_{n})\,\,\textrm{where}\,\,\theta_{k}=\sqrt{L_{k}\left[\sin^{2}\left(\frac{t-x}{2}\right),x\right]}.
\end{eqnarray*}
Then for all $f\in C_{2\pi}(\mathbb{R})$, we get
\begin{eqnarray*}
\ \parallel L_{k}(f,x)-f(x)\parallel _{2\pi}=d(\widehat{F})-o(z_{n})
\end{eqnarray*}
where $z_{n}= \max\{u_{n},v_{n}\}$.
\end{thm}
\begin{proof}
Let $f\in C_{2\pi}(\mathbb{R})$ and $x\in [-\pi,\pi]$. From
(\ref{k7}) and $(\ref{r1})$, we can write
\begin{eqnarray*}
|L_{k}(f,x)-f(x)|&\leq&L_{k}(|f(t)-f(x)|;x)+|f(x)||L_{k}(1,x)-1| \\
&\leq&L_{k}\left(\frac{|x-y|}{\delta}+1;x\right)\omega(f,\delta)+|f(x)||L_{k}(1,x)-1|\\
&\leq&L_{k}\left(\frac{\pi^{2}}{\delta^{2}}\sin^{2}\left(\frac{y-x}{2}\right)+1;x\right)\omega(f,\delta)+|f(x)||L_{k}(1,x)-1|\\
&\leq&\left\{L_{k}(1,x)+\frac{\pi^{2}}{\delta^{2}}L_{k}\left(\sin^{2}\left(\frac{y-x}{2}\right);x\right)\right\}\omega(f,\delta)+|f(x)||L_{k}(1,x)-1|\\
&=&\left\{L_{k}(1,x)+\frac{\pi^{2}}{\delta^{2}}L_{k}\left(\sin^{2}\left(\frac{y-x}{2}\right);x\right)\right\}\omega(f,\delta)+|f(x)||L_{k}(1,x)-1|.
\end{eqnarray*}
By choosing $\sqrt{\theta_{k}}=\delta$, we get
\begin{eqnarray*}
\parallel
L_{k}(f,x)-f(x)\parallel _{2\pi}&\leq& \parallel
f\parallel_{2\pi}\parallel L_{k}(1,x)-x \parallel
_{2\pi}+2\omega(f,\theta_{k})+\omega(f,\theta_{k})\parallel
L_{k}(1,x)-x
\parallel _{2\pi}\\
&\leq&K\{\parallel L_{k}(1,x)-x \parallel
_{2\pi}+\omega(f,\theta_{k})+\omega(f,\theta_{k})\parallel
L_{k}(1,x)-x
\parallel _{2\pi}\},
\end{eqnarray*}
where $K=\max\{2,\parallel f\parallel_{2\pi}\}$. By Definition
\ref{defn} and conditions (i) and (ii), we get the desired the
result.
\end{proof}

\section*{Conflict of Interests}
The authors declare that there are no conflict of interests regarding the publication of this paper.

\end{document}